\documentclass[a4paper,reqno,final]{amsart}
\usepackage[all]{xy}
\usepackage{amssymb,latexsym,amsmath}
\usepackage{amsthm}
\usepackage[colorlinks=true,breaklinks=true,linkcolor=black,citecolor=black,hidelinks]{hyperref}
\usepackage{verbatim}
\usepackage[textsize=footnotesize]{todonotes}
\usepackage{times}
\usepackage{manfnt}
\usepackage{mathrsfs}
\swapnumbers
\theoremstyle{plain}
\newtheorem{theorem}{Theorem}[section]

\newtheorem{lemma}[theorem]{Lemma}
\newtheorem{prop}[theorem]{Proposition}
\newtheorem{corollary}[theorem]{Corollary}

\theoremstyle{definition}
\newtheorem{ex}[theorem]{Example}
\newtheorem{definition}[theorem]{Definition}
\newtheorem{exer}{Exercise}
\newtheorem{rem}[theorem]{Remark}
\newtheorem{observation}[theorem]{Observation}
\newtheorem{app}[theorem]{Application}
\newtheorem{his}[theorem]{Historical Note}
\newtheorem*{notation}{Notation}
\newtheorem*{conv}{Convention}
\newtheorem{Observation}[theorem]{Observation}
\newtheorem{cons}[theorem]{Construction}

\newcommand{\bthm}{\begin{theorem}}
\newcommand{\ethm}{\end{theorem}}

\newcommand{\bprop}{\begin{prop}}
\newcommand{\eprop}{\end{prop}}

\newcommand{\bcoro}{\begin{corollary}}
\newcommand{\ecoro}{\end{corollary}}

\newcommand{\bex}{\begin{ex}}
\newcommand{\eex}{ \end{ex}}

\newcommand{\bdf}{\begin{definition}}
\newcommand{\edf}{ \end{definition}}

\newcommand{\brem}{\begin{rem}}
\newcommand{\erem}{ \end{rem}}

\newcommand{\blem}{\begin{lemma}}
\newcommand{\elem}{\end{lemma}}

\newcommand{\bhis}{\begin{his}}
\newcommand{\ehis}{ \end{his}}

\newcommand{\bobs}{\begin{observation}}
\newcommand{\eobs}{ \end{observation}}

\newcommand{\bproof}{\begin{proof}}
\newcommand{\eproof}{\end{proof}}

\newcommand{\bexr}{\begin{exer}}
\newcommand{\eexr}{ \end{exer}}

\newcommand{\bnotation}{\begin{notation}}
\newcommand{\enotation}{\end{notation}}

\newcommand{\bconv}{\begin{conv}}
\newcommand{\econv}{ \end{conv}}

\newcommand{\bapp}{\begin{app}}
\newcommand{\eapp}{ \end{app}}

\newcommand{\bfact}{\begin{Observation}}
\newcommand{\efact}{ \end{Observation}}

\newcommand{\bcons}{\begin{cons}}
\newcommand{\econs}{ \end{cons}}

\renewcommand{\bar}{\overline}

\newcommand{\Mod}{\mathsf{Mod}}
\renewcommand{\mod}{\mathsf{mod}}
\newcommand{\D}{\mathsf{D}}
\newcommand{\Add}{\mathsf{Add}}
\newcommand{\add}{\mathsf{add}}
\newcommand{\pd}{\operatorname{pd}}

\newcommand{\fd}{\operatorname{fd}}
\newcommand{\gldim}{\operatorname{gldim}}
\newcommand{\wgldim}{\operatorname{w.gldim}}
\renewcommand{\ker}{\operatorname{Ker}}
\newcommand{\coker}{\operatorname{Coker}}
\renewcommand{\hom}{\operatorname{Hom}}
\newcommand{\Ext}{\operatorname{Ext}}
\newcommand{\Tor}{\operatorname{Tor}}
\newcommand{\Rhom}{\operatorname{\mathbf{R}Hom}}
\newcommand{\tensorL}{\mathop{\otimes^{\mathbf{L}}} \limits}
\newcommand{\End}{\operatorname{End}}
\newcommand{\op}{{\operatorname{op}}}
\newcommand{\arrow}{\longrightarrow}
\renewcommand{\mapsto}{\longmapsto}
\newcommand{\im}{\operatorname{Im}}
\newcommand{\spec}{\operatorname{Spec}}
\title{\textbf{Tilting Modules under Special Base Changes}}
\author{Pooyan Moradifar \and Shahab Rajabi \and Siamak Yassemi}
\address[Pooyan Moradifar]{School of Mathematics, Statistics and Computer Science \\ College of Science \\ University of Tehran \\ Tehran \\ Iran}
\email{pmoradifar@ut.ac.ir}
\address[Shahab Rajabi]{School of Mathematics, Statistics and Computer Science \\ College of Science \\ University of Tehran \\ Tehran \\ Iran}
\email{shahabrjb@gmail.com}
\address[Siamak Yassemi]{School of Mathematics, Statistics and Computer Science \\ College of Science \\ University of Tehran \\ Tehran \\ Iran}
\email[corresponding author]{yassemi@ut.ac.ir}
\begin{document}
\begin{abstract}
Given a non-unit, non-zero-divisor, central element $x$ of a ring $\Lambda$, it is well known that many properties or invariants of $\Lambda$ determine,
and are determined by, those of $\Lambda / x \Lambda$ and $\Lambda_x$. In
the present paper, we investigate how the property of ``being tilting'' 
behaves in this situation. It turns out that any tilting module over $\Lambda$ gives rise to tilting modules over $\Lambda_x$ and $\Lambda / x \Lambda$ 
after localization and passing to quotient respectively. On the other hand, it is proved that under some mild conditions, a module over $\Lambda$ is tilting if its 
corresponding localization and quotient are tilting over $\Lambda_x$ and $\Lambda / x \Lambda$ respectively.
\end{abstract}
\maketitle
\section*{Introduction}

\bconv
Throughout this note, ``ring'' means an \emph{arbitrary ring with a non-zero unit}, unless otherwise stated. If $\Lambda$ is a ring, by a ``$\Lambda$-module'' we mean a \emph{left} $\Lambda$-module and right $\Lambda$-modules are considered as left modules over $\Lambda^\op$---the opposite ring of $\Lambda$. An element $x \in \Lambda$ is said to be \emph{regular} on a $\Lambda$-module $M$ if it is non-zero-divisor over $M$ and $xM \ne M$.
\econv

\bnotation
Given a ring $\Lambda$, the flat and the projective dimension of $\Lambda$-modules will be denoted by $\fd_\Lambda (-)$ and $\pd_\Lambda (-)$ respectively. The \emph{global dimension} and \emph{weak global dimension} of $\Lambda$ will also be denoted by $\gldim (\Lambda)$ and $\wgldim (\Lambda)$ respectively. For a central element $x \in \Lambda$, we denote by $M_x$ the localization of a $\Lambda$-module $M$ with respect to the multiplicatively closed subset $\{ x^n : n \in \mathbb{N}_0 \}$ of $\Lambda$. If $R$ is the center of $\Lambda$, then for all $\Lambda$-modules $M$ and $N$, $\Ext^n_\Lambda (M,N)$ has an $R$-module structure and so we may form its localization as an $R$-module with respect to $\{ x^n : n \in \mathbb{N}_0 \}$, denoted by $\Ext^n_\Lambda (M,N)_x$.
\enotation

Tilting theory is one of the major branches of representation theory of algebras which originally was  introduced in the context of finitely generated modules over artin algebras, mainly through the work of Brenner and Butler~\cite{brenner.butler},
Happel and Ringel~\cite{tiltedalg.happel} and Happel~\cite{happel.book}. The theory was later generalized to the setting of finitely generated 
modules over associated rings by Miyashita~\cite{miyashita}, and extended to the setting of arbitrary modules over associative rings 
by Angeleri-H\"{u}gel et al.~\cite{tilting.preenv,hugel.coel}.
One of the main themes of tilting theory is to compare the properties of a ring and its tilted ring:
Once a tilting module $T$ is detected over a ring $\Lambda$, one can look at the ``tilted ring'' $\End_{\Lambda}(T)^{\op}$ which has in general many representation-theoretic, (co)homological and K-theoretic
properties and invariants in common with the original ring $\Lambda$ (see e.g.~\cite[Section~9]{rickard.der.morita}).
This close connection between $\Lambda$ and its tilted ring is in large part due to some equivalences induced by $T$ between certain subcategories of the module categories or derived categories of the two rings; see e.g.~\cite{miyashita},~\cite{happel.onderived},~\cite{rickard.der.morita},
~\cite{bazzoni.inf.tilt1},~\cite{bazzoni.inf.tilt2} and also Remark~\ref{rem-der-equiv}.

In this note, we aim to study tilting modules under special base changes associated with a central element of a ring. Generally speaking, base change theorems deal with ascent or descent of invariants or properties of rings and their modules along a ring homomorphism $\Lambda \arrow \Lambda'$. As far as tilting theory is concerned, a basic question in this regard is that when for a 
given tilting $\Lambda$-module $T$, the $\Lambda'$-module $\Lambda' \otimes_\Lambda T$ is tilting? This problem has already been addressed---for finitely generated tilting modules---by several authors (e.g.~\cite{tachikawa.tilting},~\cite{miyachi.tilting.trivial},~\cite{miyashita},~\cite{fuller.extension} and~\cite{assem.marmaridis})
while the reverse problem, namely when  $\Lambda' \otimes_\Lambda T$ being tilting implies $T$ is tilting, has been less explored; see~\cite{trlifaj.coloc} for a relevant work in this regard.
This is partially due to the fact that one usually loses too much information (e.g. finiteness of the projective dimension) in the ``descent process''. The situation may get better if we have another ring homomorphism $\Lambda \arrow \Lambda''$ at our disposal which in a sense makes up for the ``lost information''. This situation occurs particularly for some special base changes associated with a central element of a ring: Given a central element $x$ of a ring $\Lambda$, we may form two new rings, namely $\Lambda_x$ and $\Lambda / x \Lambda$, and then we may
consider base changes along the canonical ring homomorphisms $\Lambda \arrow \Lambda_x$ and $\Lambda \arrow \Lambda / x \Lambda$. When $\Lambda$ is commutative, this corresponds---from the algebro-geometric point of view---to stratification of the affine scheme $\spec(\Lambda)$ into the open subscheme $\spec(\Lambda_x)$, which corresponds to $\{\mathfrak{p}\in \spec(\Lambda) | x \not \in \mathfrak{p}\}$, and its complement $\spec(\Lambda/x\Lambda)$, which corresponds to $\{\mathfrak{p}\in \spec(\Lambda) | x \in \mathfrak{p}\}$. The results of \cite{rajabiI} and its sequel \cite{rajabiII} bear witness to the fact that when $x$ is a non-unit, non-zero-divisor central element of a (not necessarily commutative) ring $\Lambda$, there is a close relation between homological behavior of the ring $\Lambda$ on one hand, and homological behavior of the rings $\Lambda_x$ and $\Lambda/x\Lambda$ on the other hand. The following result from~\cite{rajabiI} is an instance of this phenomenon, which will also be used in the sequel.

\bthm[{\cite[Theorem~3.2 and 3.4]{rajabiI}}] \label{thm:rajabi}
Let $M$ be a module over a ring $\Lambda$. If $x \in \Lambda$ is a non-unit central element which is non-zero-divisor on $\Lambda$ and $M$, then
\begin{enumerate}
\item $\fd_\Lambda (M) = \max \{ \fd_{\Lambda/x\Lambda} (M/xM) , \fd_{\Lambda_x} (M_x)  \}$. If furthermore $M$ has a degreewise finitely generated projective resolution,
then we may replace `flat dimension'' with ``projective dimension'' in the equality. 

\item $\wgldim (\Lambda) \leq \max \big \{ 1+\wgldim (\Lambda / x \Lambda) , \wgldim (\Lambda_x) \big \}$.

\item $\gldim (\Lambda) \leq \max \big \{ 1+\gldim (\Lambda/x\Lambda) , \gldim (\Lambda_x) \big  \}$ provided that $\Lambda$ is left noetherian.
\end{enumerate}
Furthermore, the equality in parts~(ii) and~(iii) is attained provided that $\wgldim (\Lambda / x \Lambda)$ and $\gldim (\Lambda/x\Lambda)$ are finite respectively. 
\ethm

 A tilting module is a module of finite projective dimension with some additional properties (see~\ref{df:tilting}), and we know
in light of Theorem~\ref{thm:rajabi} that under certain mild conditions, finiteness of the projective dimension 
of a module $T$ over $\Lambda$ forces and is forced by finiteness 
of the projective dimension of the ``extended'' modules $T_x$ and $T/xT$ over $\Lambda_x$ and $\Lambda / x \Lambda$ respectively. This is our point of departure in studying tilting modules under special
base changes along $\Lambda \arrow \Lambda_x$ and $\Lambda \arrow \Lambda /x \Lambda$.
Our results can be divided into two parts depending on the type of the tilting module under consideration, i.e. ``classical'' or ``non-classical''. The classical case is discussed in Theorem~\ref{main.thm.classic}, where we prove that under mild conditions the following bi-implication holds for a $\Lambda$-module $T$:
\[\text{$T$ is a classical $n$-tilting $\Lambda$-module} \Longleftrightarrow \left \{ 
\begin{array}{l}
\text{$T_x$ is a classical $n$-tilting $\Lambda_x$-module} \\[0.3cm]
\text{$T/xT$ is a classical $n$-tilting $\Lambda/x \Lambda$-module}
\end{array} \right. \]
Furthermore, it is proved that $\End_{\Lambda_x} (T_x)^\op \cong \Gamma_{x_\bullet}$ and $\End_{\Lambda /  x \Lambda} (T/xT)^\op \cong \Gamma / x_\bullet \Gamma$, where $\Gamma:=\End_\Lambda (T)^\op$ is the tilted ring of $\Lambda$ via $T$, and $x_\bullet$ is the central element of $\Gamma$ defined by multiplication by $x$. It is well known that classical tilting modules induce equivalences on the level of derived categories (see e.g.~\cite{happel.book} or~\cite{keller.Handbook}). In particular, when the above bi-implication holds, there are derived equivalences 
$\D (\Lambda) \simeq \D (\Gamma)$, $\D (\Lambda_x) \simeq \D (\Gamma_{x_\bullet})$
and $\D (\Lambda/x \Lambda) \simeq \D (\Gamma/x_{\bullet} \Gamma)$, which implies a strong connection between $\Lambda$ and $\Gamma$ 
as well as their corresponding localizations and quotients; see Remark~\ref{rem-der-equiv} for more details.

Regarding the non-classical case, the problem and of course its required arguments are more subtle: It is proved in Theorem~\ref{main.thm.nonclassic.loc}
that when $x$ is regular on a $\Lambda$-module $T$, the implication
\[ \text{$T$ is an $n$-tilting module} \implies \text{$T_x$ and $T/xT$ are $n$-tilting modules} \]  
holds. The reverse implication, however, is \emph{not} true in general. In fact, as Example~\ref{counterexample} shows, even the following weaker form of the reverse implication does \emph{not} hold in general.
\[ \text{$ T_x$ and $T/xT$ are tilting modules} \implies \text{$T$ is a partial tilting module}.  \]  
In view of Theorem~\ref{main.thm:tilting.class}, this failure is in large part
due to the fact that the class $T^{\perp_\infty}$ might fail to be closed under direct sums. Indeed, we prove in Theorem~\ref{main.thm:tilting.class} that
 when $T^{\perp_\infty}$ is closed under direct sums, the above-mentioned implication holds under some reasonable conditions. In particular, $T$ is a direct summand of a tilting module.

As a final introductory remark, it is worth mentioning that the relation between $\Lambda$, $\Lambda_x$ and $\Lambda / x \Lambda$, where $x$ is a regular central element of $\Lambda$, is in many respects analogous to the relation between the rings $\Lambda$, $x\Lambda x$ and $\Lambda / \Lambda x \Lambda$ where $x$ is an idempotent element. Indeed, in this case we can still form the ring $\Lambda_x$ and it is easily seen that $\Lambda_x \cong x\Lambda x$, and the ring $\Lambda / x \Lambda$ corresponds to the ring $\Lambda / \Lambda x \Lambda$.  The interplay between various invariants and properties of $\Lambda$, $x \Lambda x$ and $\Lambda / \Lambda x\Lambda$, where $x$ is an idempotent element, is extensively explored by Cline, Parshall and Scott
via recollements of derived categories; see e.g.~\cite{CP.quasi}, \cite{CPS.highest} and~\cite{CPS}. In fact, there is a necessary and sufficient 
condition for existence of a recollement of $\D_\ast (\Lambda)$ relative to $\D_\ast (\Lambda/\Lambda x \Lambda)$ and $\D_\ast (x \Lambda x)$ when $x$ is an idempotent element; see~\cite{CPS}.
However, to the best of our knowledge, the existence of a similar recollement in our setting, where $x$ is a regular central element, is unknown in general.


\section{Preliminaries}

In this section, we recall some notions and facts from approximation theory and tilting theory, which will be used in the sequel.
We begin with some notational remarks about syzygies and orthogonal classes.

\bnotation
Given a ring $\Lambda$, the category of $\Lambda$-modules is denoted by $\Mod (\Lambda)$, and the full subcategory of $\Mod (\Lambda)$ consisting 
of \emph{strongly finitely presented modules} (i.e. $\Lambda$-modules possessing a degreewise finitely generated projective resolution)
is denoted by $\mod (\Lambda)$. Furthermore, for every integer $n \geq 0$, the category of $\Lambda$-modules of projective dimension at most $n$ is 
denoted by $\mathcal{P}_n$ and we also let $\mathcal{P}_n^{< \omega} := \mathcal{P}_n \cap \mod (\Lambda)$. 

Let $M$ be a $\Lambda$-module and 
\[\xymatrix{ \mathbf{P} = \cdots \ar[r] & P_1 \ar[r]^-{d_1} & P_0 \ar[r]^-{d_0} & M \ar[r] & o  }\]
be an augmented projective resolution of $M$. For every integer $n \geq 0$, the $\Lambda$-module $\Omega_i (\mathbf{P}) : = \im d_n$ is called the
\emph{$n$-th syzygy} module of $M$ with respect to the projective resolution $\mathbf{P}$.

Let $\Lambda$ be a ring. For a class $\mathcal{X}$ of $\Lambda$-modules let
\[ \mathcal{X}^{\perp} : = \big \{ M \in \Mod (\Lambda) : \text{$\Ext^1_\Lambda (X , M) = o$ for all $X \in \mathcal{X}$} \big \} \: , \]
and
\[ {}^{\perp} \mathcal{X} : = \big \{ M \in \Mod (\Lambda) : \text{$\Ext^1_\Lambda (M,X) = o$ for all $X \in \mathcal{X}$} \big \} \: . \]
These are called the \emph{right} and the \emph{left} \emph{orthogonal classes} of $\mathcal{X}$ respectively.
In the special case where $\mathcal{X}$ consists of only one $\Lambda$-module $T$, we let
$T^\perp := \{ T \}^\perp$ and ${}^\perp T :={}^\perp \{ T \}$ for convenience.
\enotation

\bdf \label{df:filtration}
Let $\Lambda$ be a ring. A \emph{continuous chain} of $\Lambda$-modules (indexed by an ordinal $\kappa$) is a family $\{ M_\alpha \}_{\alpha \leq \kappa}$ of $\Lambda$-modules such that $M_0 = o$, $M_{\alpha} \subseteq M_{\alpha +1}$ for any 
ordinal $\alpha < \kappa$, and $M_\alpha = \bigcup_{\beta < \alpha} M_\beta$ for any limit ordinal $\alpha \leq \kappa$. 

Given a class $\mathcal{C}$ of $\Lambda$-modules, a $\Lambda$-module $M$ is said to be \emph{$\mathcal{C}$-filtered} if there is a continuous chain $\{ M_{\alpha} \}_{\alpha \leq \kappa}$ of submodules of $M$ such that $M_{\kappa} = M$ and for every ordinal $\alpha < \kappa$, the $\Lambda$-module $M_{\alpha +1} / M_{\alpha}$ is isomorphic to some element of $\mathcal{C}$. The family $\{ M_{\alpha} \}_{\alpha \leq \kappa}$ is then called a \emph{$\mathcal{C}$-filtration} of $M$.
\edf

The most important closure property of left orthogonal classes is that they are closed under filtration, which is the theme of the following well-known lemma. 

\blem[{Eklof~\cite{eklof}}] \label{Eklof.lemma}
Let $\Lambda$ be a ring and $N$ be a $\Lambda$-module. If a $\Lambda$-module $M$ is ${}^\perp N$-filtered, 
then $M \in {}^\perp N$.
\elem

\bproof
See~\cite[Lemma 6.2]{trlifaj}.
\eproof

Approximation theory of modules was founded independently by Auslander et al.~\cite{aus.preproj,AR.applications} and Enochs~\cite{enochs.coversenv}. The theory has many applications in representation theory of algebras in general, and tilting theory in particular; see~\cite[Part~(III)]{trlifaj} for more information in this regard. Here we recall some notions and facts from approximation theory which will be needed in the sequel.

\bdf \label{df:special.precov} \label{df:special.preenv} 
Let $\Lambda$ be a ring and $\mathcal{X}$ be a class of $\Lambda$-modules. A $\Lambda$-monomorphism $f : M \arrow X$ with $X \in \mathcal{X}$ is called a \emph{special $\mathcal{X}$-preenvelope} of $M$ if $\coker f \in {}^\perp \mathcal{X}$. 
The class $\mathcal{X}$ is said to be \emph{special preenveloping}
if every $\Lambda$-module has a special $\mathcal{X}$-preenvelope. The notion of a ``special $\mathcal{X}$-precover'' and a ``special precovering class'' is defined dually; cf.~\cite[Chapter~5]{trlifaj}.
\edf

The following celebrated theorem of Eklof and Trlifaj~\cite{trlifaj.Extvanish} states that special preenvelopes are abundant.

\bthm \label{thm:eklof.trlifaj}
Let $\Lambda$ be a ring, $\mathcal{S}$ be a \emph{set} of $\Lambda$-modules, and $\mathcal{S}^+:=\mathcal{S} \cup \{ \Lambda \}$.
\begin{enumerate}
\item For every $\Lambda$-module $M$ there is a short exact sequence
\[\xymatrix{ \delta= o \ar[r] & M \ar[r]^-f & P  \ar[r] & N \ar[r] & o }\]
of $\Lambda$-homomorphism wherein $P \in \mathcal{S}^\perp$ and $N$ is $\mathcal{S}$-filtered. In particular, $f : M \arrow P$ is a 
special $\mathcal{S}^\perp$-preenvelope of $M$.

\item The class ${}^\perp (\mathcal{S}^\perp)$ consists precisely of all $\Lambda$-modules which are direct summands of 
$\mathcal{S}^+$-filtered modules.
\end{enumerate}
\ethm

\bproof
See~\cite[Theorem 6.11]{trlifaj} and~\cite[Corollary 6.13]{trlifaj}.
\eproof

\bnotation \label{not:inf.classes}
Let $\Lambda$ be a ring. For a class $\mathcal{X}$ of $\Lambda$-modules let
\begin{align*}
\mathcal{X}^{\perp_\infty} & : =  \big \{ M \in \Mod (\Lambda) : \hbox{$\Ext^{\geqslant 1}_\Lambda (X , M) = o$ for all $X \in \mathcal{X}$}   \big \} \: , \\
{}^{\perp_\infty} \mathcal{X} & : =  \big \{ M \in \Mod (\Lambda) : \hbox{$\Ext^{\geqslant 1}_\Lambda (M,X) = o$ for all $X \in \mathcal{X}$}   \big \} \: .
\end{align*}
We also let $T^{\perp_\infty} := \{ T \}^{\perp_\infty}$ and ${}^{\perp_\infty} T :={}^{\perp_\infty} \{ T \}$ for convenience.
\enotation

\bcoro \label{coro:T.perp.preenv}
Let $T$ be a module over a ring $\Lambda$. If $\mathbf{\Omega}$ is the set of all syzygy modules of $T$ with respect to some projective resolution $\mathbf{P}$ of $T$, then $T^{\perp_\infty}=\mathbf{\Omega}^\perp$. In particular, $T^{\perp_\infty}$ is special preenveloping.
\ecoro

\bproof
The equality $T^{\perp_\infty} =\mathbf{\Omega}^\perp$ holds by the fact that
\[\Ext^{n+1}_\Lambda \big (T , M \big) \cong \Ext^1_\Lambda \big ( \Omega_{n} (\mathbf{P}) , M \big )\]
for every integer $n \geq 0$. Since $\mathbf{\Omega}$ is a set, it follows from Theorem~\ref{thm:eklof.trlifaj} that $T^{\perp_\infty}$ is special preenveloping.
\eproof


\bdf \label{df:tilting}
Let $\Lambda$ be a ring and $n \geq 0$ be an integer. A $\Lambda$-module $T$ is called an \emph{$n$-tilting module} if
\begin{description}
\item[(T1)] $\pd_\Lambda (T) \leq n$.

\item[(T2)] For every cardinal number $\kappa$, $\Ext^{\geqslant 1}_\Lambda \big (T,T^{(\kappa)} \big )=o$ (here $T^{(\kappa)}$ denotes the direct sum of $\kappa$ copies of $T$).

\item[(T3)] There is an exact sequence
\[\xymatrix{o \ar[r] & \Lambda \ar[r] & T_0 \ar[r] & \cdots \ar[r] & T_m \ar[r] & o}\]

wherein $T_i \in \Add_\Lambda (T)$ for all $0 \leq i \leq m$.
\end{description}
Here $\Add_\Lambda (T)$ denotes the class of all $\Lambda$-modules which are direct summands of direct sums of copies of $T$. 
When $T$ is an $n$-tilting module,  the ring $\Gamma:=\End_\Lambda (T)^\op$ is called the \emph{tilted ring} of $\Lambda$ via $T$, and
the triple $(\Lambda , T , \Gamma)$ is called the \emph{$n$-tilting triple} associated with $T$. Note that $T$ is a ($\Lambda$-$\Gamma$)-bimodule in the natural way. A $\Lambda$-module is called \emph{tilting} if it is $n$-tilting for some integer $n \geq 0$.

An $n$-tilting $\Lambda$-module $T$ is called a \emph{classical $n$-tilting module} if it is strongly finitely presented. It can be proved (see e.g.~\cite[page 301]{trlifaj})
that an $n$-tilting $\Lambda$-module $T$ is classical $n$-tilting if and only if it satisfies the following ``weaker'' form of (T1)--(T3):
\begin{description}
\item[(t1)] $T$ possesses a degreewise finitely generated projective resolution of length at most $n$.

\item[(t2)] $\Ext^{\geqslant 1}_\Lambda (T,T)=o$.

\item[(t3)] There exists an exact sequence
\[\xymatrix{o \ar[r] & \Lambda \ar[r] & T_0 \ar[r] & \cdots \ar[r] & T_m \ar[r] & o}\]
of $\Lambda$-modules where $T_i \in \add_\Lambda (T)$ for  $i=0, \ldots , m$.
\end{description}
Here $\add_\Lambda (T)$ denotes the class of all $\Lambda$-modules which are direct summands of finite direct sums of copies of $T$. 
The tilting triple associated with a classical $n$-tilting module is called a \emph{classical $n$-tilting triple}.
\edf

The following characterization of classical tilting modules, implicit in the work of Miyashita~\cite[Section~1]{miyashita} and easily deduced from
Wakamatsu~\cite[Proposition~5]{wakamatsu}, will be useful in the sequel.

\bprop \label{prop:tilting.characterization}
Let $\Lambda$ be a ring and $T$ be a $\Lambda$-module. 
Let $\Gamma:=\End_\Lambda (T)^\op$ and regard $T$ as a ($\Lambda$-$\Gamma$)-bimodule in the natural way. Then the $\Lambda$-module $T$ is a classical $n$-tilting module if and only if the following conditions are satisfied:
\begin{enumerate}
\item $T$ has a degreewise finitely generated projective resolution of length at most $n$ as a $\Lambda$-module and as a $\Gamma^\op$-module.

\item $\Ext^{\geqslant 1}_{\Lambda} (T,T)=\Ext^{\geqslant 1}_{\Gamma^\op} (T,T)=o$.

\item The natural ring homomorphism $\Lambda \arrow \End_{\Gamma^\op} (T)$ is an isomorphism.
\end{enumerate}
\eprop

As we explained in the introduction, in this note we are mainly concerned with a particular case of the following general ``base change problem'': Given a ring homomorphism
$\Lambda \arrow \Lambda'$, under what conditions a tilting $\Lambda$-module $T$ induces an ``extended'' tilting $\Lambda'$-module $\Lambda' \otimes_\Lambda T$? 
To the best of our knowledge, the most general result in this regard is due to Miyashita~\cite{miyashita} as follows.

\bprop \label{prop:miyashita.base.change}
Let $\Lambda \arrow \Lambda'$ be a ring homomorphism and assume that $T$ is a $\Lambda$-module such 
that $\Tor^{\Lambda}_{\geqslant 1} (\Lambda' , T) = o$. 
\begin{enumerate}
\item  If $T$ is a classical $n$-tilting module and $\Ext^{\geqslant 1}_{\Lambda'} (\Lambda' \otimes_\Lambda T , \Lambda' \otimes_\Lambda T )=o$,
then $\Lambda' \otimes_\Lambda T$ is a classical $n$-tilting $\Lambda'$-module.

\item  If $T$ is an $n$-tilting module 
and $\Ext^{\geqslant 1}_{\Lambda'} (\Lambda' \otimes_\Lambda T , \Lambda' \otimes_\Lambda T^{(\kappa)} )=o$ for any cardinal number $\kappa$,
then $\Lambda' \otimes_\Lambda T$ is an $n$-tilting $\Lambda'$-module.
\end{enumerate}
\eprop

\bproof
The result is proved for classical tilting modules in~\cite[Theorem~5.2]{miyashita}, and the same argument also applies to non-classical tilting modules.
\eproof

\bdf \label{df:tilting.class}
Let $\Lambda$ be a ring. A class $\mathcal{C}$ of $\Lambda$-modules is called an \emph{$n$-tilting class} (for some integer $n \geq 0$)
if there exists an $n$-tilting module $T$ with $\mathcal{C} = T^{\perp_\infty}$. 
The class $\mathcal{C}$ is called \emph{tilting} if it is $n$-tilting for some integer $n \geq 0$.
\edf

It is easily seen that every tilting class is \emph{coresolving} in the sense that it is closed under extensions and cokernel of monomorphisms, 
and it contains all injective $\Lambda$-modules. The following result, due to Angeleri-H\"{u}gel and Coelho~\cite{hugel.coel}, characterizes tilting classes among
coresolving classes in $\Mod (\Lambda)$.

\bthm \label{thm:tilting.class}
Let $\Lambda$ be a ring and $n \geq 0$ be an integer. A class $\mathcal{C}$ of $\Lambda$-modules is $n$-tilting if and only if the following conditions are satisfied:
\begin{enumerate}
\item $\mathcal{C}$ is coresolving in $\Mod (\Lambda)$, closed under direct sums and direct summands,

\item ${}^\perp \mathcal{C} \subseteq \mathcal{P}_n$,

\item $\mathcal{C}$ is special preenveloping.
\end{enumerate}
\ethm

\bproof
See~\cite[Theorem 13.18]{trlifaj}.
\eproof

\bdf \label{df:partial.tilting}
Let $\Lambda$ be a ring and $n \geq 0$ be an integer. A $\Lambda$-module $P$ is called a \emph{partial $n$-tilting module} provided that
\begin{description}
\item[(P1)] $\pd_\Lambda (P) \leq n$,

\item[(P2)] $\Ext^{\geqslant 1} (P,P^{(\kappa)})=o$ for any cardinal number $\kappa$.

\item[(P3)] $P^{\perp_\infty}$ is closed under direct sums.
\end{description}
\edf

\brem \label{rem:partial}
Given a partial tilting $\Lambda$-module $P$,
it follows directly from Theorem~\ref{thm:eklof.trlifaj} and Theorem~\ref{thm:tilting.class} that $P^{\perp_\infty}=T^{\perp_\infty}$ for some $n$-tilting module $T$. With a little work one can deduce from Theorem~\ref{thm:tilting.class} that $T$ can be chosen in such a way that
it contains $P$ as a direct summand; see~\cite[Corollary 13.22]{trlifaj}.
\erem

One of the most significant properties of tilting classes is that they are of ``finite type'' (see~\cite{tilting.ftype}) in the sense that they are of the form $\mathcal{S}^{\perp_\infty}$ for some $\mathcal{S} \subseteq \mathcal{P}^{< \omega}_n$.
This result yields the following well-known correspondence between tilting classes and certain resolving subcategories of $\mod (\Lambda)$
which will be used later in proving our main results.
Recall that a class $\mathcal{S}$ of $\Lambda$-modules is called a \emph{resolving subcategory of $\mod (\Lambda)$} if $\mathcal{P}^{< \omega}_0 \subseteq \mathcal{S} \subseteq \mod (\Lambda)$, $\mathcal{S}$ is closed under direct sums, direct summands and kernel of monomorphisms.

\bthm \label{thm:finite.type}
Let $\Lambda$ be a ring and $n \geq 0$ be an integer. Then there is a one-to-one correspondence between
$n$-tilting classes $\mathcal{T}$ and resolving subcategories $\mathcal{S}$ of $\mod (\Lambda)$ such that $\mathcal{S} \subseteq \mathcal{P}^{< \omega}_n$. The correspondence is given by the mutually inverse assignments
$\mathcal{T} \mapsto ({}^\perp \mathcal{T}) \cap \mod (\Lambda)$ and $\mathcal{S} \mapsto \mathcal{S}^\perp$.
\ethm

\bproof
See~\cite[Theorem 13.49]{trlifaj}.
\eproof

\section{Main Results}

Let $\Lambda$ be a ring and $x \in \Lambda$ be a non-unit central element. In this section we focus on the subject matter of this paper, 
namely studying tilting modules under special base changes along the canonical ring homomorphisms $\Lambda \arrow \Lambda_x$ and
$\Lambda \arrow \Lambda / x \Lambda$. Our results can be divided into two parts, depending on the type of the tilting module under consideration, i.e. classical or non-classical. Our first main result concerning classical tilting modules under special bases changes is as follows.

\begin{theorem} \label{main.thm.classic}
Let $\Lambda$ be a ring and fix a central element $x \in \Lambda$. Let $T$ be a strongly finitely presented $\Lambda$-module over which $x$ is regular, and let $\Gamma:=\End_\Lambda (T)^\op$. 
Denote by $x_\bullet$ the central element of $\Gamma$ defined by $t \mapsto x t$ for every $t \in T$.
\begin{enumerate}
\item \sloppy If $(\Lambda , T , \Gamma)$ is a classical $n$-tilting triple, then $\big (\Lambda/x\Lambda , T/xT , \Gamma / x_\bullet \Gamma \big )$
and $(\Lambda_x , T_x , \Gamma_{x_\bullet})$ are classical $n$-tilting triples.

\item Assume that $\Lambda$ is finitely generated as a module over a noetherian center and $x$ is also regular on $\Lambda$.
If the $\Lambda_x$-module $T_x$ and the $\Lambda / x \Lambda$-module $T/xT$ are classical $n$-tilting, then $T$ is a classical $n$-tilting $\Lambda$-module.
\end{enumerate}
\end{theorem}

Before we proceed to prove Theorem~\ref{main.thm.classic}, we need to prove some preparatory results. We begin with the following two lemmas which will help us to relate vanishing of $\Ext$-modules over $\Lambda$ to vanishing of $\Ext$-modules 
over $\Lambda_x$ and $\Lambda/ x \Lambda$.

\blem \label{lem:Ext.reg.element} 
Let $M$ be a module over a ring $\Lambda$ and denote the center of $\Lambda$ by $R$. Let $x \in R$ be regular on both $\Lambda$ and $M$, and let $x_\bullet : M \arrow M$ be the $\Lambda$-homomorphism defined by multiplication by $x$.
\begin{enumerate}
\item If $N$ is a $\Lambda$-module with $x N=o$, then $\Ext^n_\Lambda (M,N) \cong_R \Ext^n_{\bar{\Lambda}} (M/xM , N)$ for all $n \geq 0$.

\item If $\Ext^1_\Lambda (M,M)=o$, then $\frac{\End_\Lambda (M)}{x_\bullet \End_{\Lambda} (M)} \cong \End_{\bar{\Lambda}} (M/xM)$ canonically  as rings and $R/xR$-modules.
\end{enumerate}
\elem

\bproof
Part~(i): Let $\bar{\Lambda}=\Lambda / x \Lambda$ and note that since $xN=o$, the $\Lambda$-module $N$ naturally acquires a $\bar{\Lambda}$-module structure.
Let $\mathbf{P}$ be a projective resolution of $M$. Since $x$ is regular on both $\Lambda$ and $M$, $\Tor^\Lambda_{\geqslant 1} (\bar{\Lambda} , M)=o$. Therefore $\bar{\Lambda} \otimes_\Lambda \mathbf{P}$ is a projective resolution of the $\bar{\Lambda}$-module $M/xM$.
Consequently, 
\begin{align*}
\Ext^n_{\Lambda} (M,N) & = H_{-n} \big ( \hom_\Lambda ( \mathbf{P} , N ) \big ) \\
& \cong_R H_{-n} \Big ( \hom_\Lambda \big ( \mathbf{P} , \hom_{\bar{\Lambda}} (\bar{\Lambda} , N ) \big ) \Big ) \\
& \cong_R H_{-n} \big ( \hom_{\bar{\Lambda}} (\bar{\Lambda} \otimes_\Lambda  \mathbf{P} , N  ) \big ) \\
& \cong_R \Ext^n_{\bar{\Lambda}} (M/xM,N) \: ,
\end{align*}
for all $n \geq 0$, as desired.

Part~(ii): Applying the functor $\hom_\Lambda (M,-)$ to the exact sequence
\[\xymatrix{ o \ar[r] & M \ar[r]^-x & M \ar[r] & M/xM \ar[r] & o }\]
of $\Lambda$-modules yields the exact sequence
\[\xymatrix{ o \ar[r] & \End_\Lambda (M) \ar[r]^-{x_\bullet} & \End_\Lambda (M) \ar[r] & \hom_\Lambda (M,M/xM) \ar[r] & o  }\]
of $R$-modules in view of $\Ext^1_\Lambda (M,M)=o$. Therefore, 
\begin{align*}
\End_\Lambda (M) / x_\bullet \End_\Lambda (M) & \cong_R \hom_\Lambda (M, M/xM) \\
& \cong_R \End_{\bar{\Lambda}} (M/xM) \: ,
\end{align*}
where the second isomorphism is given by $\hom$-Tensor Adjunction. It is readily seen that the canonical
 $R$-isomorphism $\End_\Lambda (M) / x_\bullet \End_\Lambda (M) \cong \End_{\bar{\Lambda}} (M/xM)$ is both a ring isomorphism and $R/xR$-isomorphism.
\eproof

\blem \label{lem:Ext.locx}
Let $\Lambda$ be a ring with the center $R$ and $x \in R$. Let $M$ and $N$ be $\Lambda$-modules and, $\bar{\Lambda}:=\Lambda/x \Lambda$.
if $M$ is strongly finitely presented, then there is a natural $R_x$-isomorphism
\[ \Ext^i_\Lambda (M,N)_x \cong \Ext^i_{\Lambda_x} (M_x,N_x) \: . \]
for every integer $i \geq 0$. The isomorphism is also a ring isomorphism when $i=0$ and $M=N$.
\elem

\bproof
Let us first prove the assertion for $i=0$. In this case we should prove that there is a natural $R_x$-isomorphism
\[ \hom_\Lambda (M,N)_x \cong \hom_{\Lambda_x} (M_x,N_x) \: . \]
Note that the well-known natural isomorphisms
\[ \hom_\Lambda (\Lambda , N)_x \cong_{\Lambda_x} N_x \cong_{\Lambda_x}  \hom_{\Lambda_x} (\Lambda_x,N_x) \]
give rise to the natural $R_x$-isomorphisms
\begin{equation} \tag{$\ast$}
\hom_\Lambda (\Lambda^n , N)_x \cong N^n_x \cong \hom_{\Lambda_x} (\Lambda^n_x,N_x)
\end{equation}
for every integer $n \geq 0$. Since $M$ strongly finitely presented, there exists an exact sequence
\[\xymatrix{ \Lambda^m \ar[r] & \Lambda^n \ar[r]  & M \ar[r] & o }\]
of $\Lambda$-modules, which yields the commutative diagram
\[\xymatrix{
o \ar[r] & \hom_\Lambda (M,N)_x \ar[r] & \hom_\Lambda (\Lambda^n , N)_x \ar[r] \ar[d]^-\cong & \hom_\Lambda (\Lambda^m , N)_x \ar[d]^-\cong  \\
o \ar[r] & \hom_{\Lambda_x} (M_x , N_x) \ar[r] & \hom_{\Lambda_x} (\Lambda^n_x , N_x) \ar[r]  & \hom_{\Lambda_x} (\Lambda^m_x , N_x) 
}\]
of $R_x$-modules with exact rows, wherein the vertical $R_x$-isomorphisms are the natural isomorphisms in~($\ast$). These $R_x$-isomorphisms induce a natural $R_x$-isomorphism
\[\theta : \hom_\Lambda (M,N)_x \stackrel{\cong}{\arrow} \hom_{\Lambda_x} (M_x , N_x) \: . \]
Every element of $\hom_\Lambda (M,N)_x$ is of the form $\frac{f}{x^k}$, where $f \in \hom_\Lambda (M,N)$ and
$k \geq 0$ is an integer, and it is straightforward to check that the $\Lambda_x$-homomorphism $\theta (f / x^k) : M_x \arrow N_x$ is defined by
$a/x^m \mapsto f(a) / x^{m+k}$. Consequently, $\theta$ is also a ring homomorphism when $M=N$.

As for the proof of the case $i > 0$, let $\mathbf{P}$ be a degreewise finitely generated projective resolution of $M$ and note that $\mathbf{P}_x$ is a projective resolution of the $\Lambda_x$-module $M_x$ by flatness of $\Lambda_x$ over $\Lambda$. Now by validity of the assertion for the case $i=0$, we have the $R_x$-isomorphisms
\begin{align*}
\Ext^i_{\Lambda_x} (M_x,N_x) & = H_{-i} \big ( \hom_{\Lambda_x} (\mathbf{P}_x , N_x) \big ) \\
& \cong   H_{-i} \big ( \hom_{\Lambda} (\mathbf{P} , N)_x \big ) \\
& \cong  \big (H_{-i} \hom_{\Lambda} (\mathbf{P} , N) \big )_x \\
& \cong  \Ext^i_\Lambda (M,N)_x  \: .
\end{align*}
This completes the proof.
\eproof

\bobs \label{obs:zerox}
Let $\Lambda$ be a ring and $x \in \Lambda$ be a central element. \emph{If $M$ is a finitely generated $\Lambda$-module, then $M=o$ if 
and only if $M_x = o$ and $M/xM=o$.} Indeed, since $M$ is finitely generated, the equality $M_x=o$ implies that there exists an integer $n >0$ such that $x^n M=o$. On the other hand, $M/xM=o$ amounts to $M=xM$. Therefore, $M=x^n M=o$. The reverse implication holds obviously.
\eobs

The following proposition explains when
self-orthogonality of a module over $\Lambda$ can be detected from self-orthogonality of its corresponding localization and
quotient over $\Lambda_x$ and $\Lambda / x \Lambda$ respectively.

\bprop \label{prop:self.orthogonal}
Let $T$ and $M$ be modules over a ring $\Lambda$.
Assume that $T$ is finitely generated as a module over its left-noetherian endomorphism ring and the $\Lambda$-module $M$ is strongly finitely generated. 
Let $x \in \Lambda$ be a central element and $\bar{\Lambda}:=\Lambda / x \Lambda$. If $x$ is non-zero-divisor over $M$, $T$ and $\Lambda$, then for every integer $n \geq 0$,
$\Ext^{n}_{\Lambda_x} (M_x , T_x) = \Ext^{n}_{\bar{\Lambda}} (M/xM , T/xT)=o$ implies $\Ext^{n}_\Lambda (M,T)=o$.
\eprop

\bproof
Let $\Gamma : = \End_\Lambda (T)$ and denote by $x_\bullet$ the central element of $\Gamma$ defined by $t \mapsto x t$ for all $t \in T$.
Notice that for every integer $n$ the abelian group $\Ext^n_\Lambda (M,T)$ has a $\Gamma$-module structure induced by the $\Gamma$-module structure of $T$, and it follows from the assumptions that the 
$\Gamma$-module $\Ext^n_\Lambda (M,T)$ is finitely generated.
Now, since by the hypothesis $x$ is a non-zero-divisor over $T$, the short exact sequence
\[\xymatrix{
o \ar[r] & T \ar[r]^-x & T \ar[r] & T / x T \ar[r] & o
}\]
yields the exact sequence
\begin{equation} \label{eq:prop.Extx.zero}
\xymatrix{
\Ext^n_\Lambda (M,T) \ar[r]^-{x_\bullet} & \Ext^n_\Lambda (M,T) \ar[r] &
\Ext^n_\Lambda (M, T / x T)}
\end{equation}
for every integer $n \geq 0$. Since $\Ext^n_\Lambda (M, T / x T)=\Ext^n_{\bar{\Lambda}} (M/xM , T/xT)=o$ by Lemma~\ref{lem:Ext.reg.element}-(i), it follows from the exact sequence~\eqref{eq:prop.Extx.zero} that \sloppy
$\Ext^n_\Lambda (M,T)=x_\bullet \Ext^n_\Lambda (M,T)$.
On the other hand, $\Ext^n_\Lambda (M,T)_{x_\bullet}=o$ by the hypothesis. Therefore, $\Ext^n_\Lambda (M,T)=o$ by Observation~\ref{obs:zerox}.
\eproof

\brem \label{rem:on.prop}
The assertion of Proposition~\ref{prop:self.orthogonal} is valid in particular when $\Lambda$ is finitely generated as a module over a notherian center, and the $\Lambda$-modules $T$ and $M$ are finitely generated.
\erem

We are now in a position to prove Theorem~\ref{main.thm.classic}.

\bproof[\textsc{Proof of Theorem~\ref{main.thm.classic}}]
Let $\bar{\Lambda} : = \Lambda / x \Lambda$ and $\bar{T}:= T / x T$. Denote the center of $\Lambda$ by $R$.

Part~(i): Since $\Ext^{1}_\Lambda (T,T)=o$, it follows from Lemma~\ref{lem:Ext.reg.element}-(ii) and~\ref{lem:Ext.locx}
that there are canonical ring isomorphisms $\Gamma / x_\bullet \Gamma \cong \End_{\bar{\Lambda}} (T/xT)^\op$ and $\Gamma_{x_\bullet} \cong \End_{\Lambda_x} (T_x)^\op$.
Notice that since $T$ is a tilting module and $x$ is regular on $T$,
the property (T3) in the definition of a tilting module implies that $x$ is regular on $\Lambda$ too. This in conjunction with flatness of $\Lambda_x$ over $\Lambda$ implies 
that $\Tor^{\Lambda}_{\geqslant 1} (\bar{\Lambda} , T)=\Tor^{\Lambda}_{\geqslant 1} (\Lambda_x , T)=o$. Therefore, by Proposition~\ref{prop:miyashita.base.change}, in order to prove part~(i) it suffices to show that $\Ext^{\geqslant 1}_{\Lambda_x} (T_x, T_x)=\Ext^{\geqslant 1}_{\bar{\Lambda}} (\bar{T} , \bar{T})=o$. Note that
\begin{align}
\Ext^n_{\Lambda_x} (T_x , T_x)  & \cong_{R_x} \Ext^n_{\Lambda} (T ,T)_x \: , \label{eq:Extx}  \\[0.3cm] 
\Ext^{n}_{\bar{\Lambda}} (\bar{T} , \bar{T}) & \cong_{R} 
\Ext^{n}_{\Lambda} (T , \bar{T}) \: , \label{eq:Extbar}
\end{align}
for every integer $n \geq 1$, where~\eqref{eq:Extx} follows from Lemma~\ref{lem:Ext.locx} and~\eqref{eq:Extbar} follows from Lemma~\ref{lem:Ext.reg.element}-(i). Now~\eqref{eq:Extx} in view of $\Ext^{\geqslant 1}_\Lambda (T,T)=o$ implies $\Ext^{\geqslant 1}_{\Lambda_x} (T_x, T_x)=o$.
On the other hand, the short exact sequence
\[\xymatrix{ o \ar[r] & T \ar[r]^-x & T \ar[r] & T/xT \ar[r] & o }\]
of $\Lambda$-modules induces the exact sequence
\begin{equation} \label{eq:Exttemp}
\xymatrix{\Ext^n_\Lambda (T,T)  \ar[r] & \Ext^n_\Lambda ( T, \bar{T} ) \ar[r] & \Ext^{n+1}_\Lambda (T,T) 
}\end{equation}
of $R$-modules for every integer $n \geq 1$. Now~\eqref{eq:Exttemp} in view of $\Ext^{\geqslant 1}_\Lambda (T,T)=o$ and~\eqref{eq:Extbar}
implies that $\Ext^{\geqslant 1}_{\bar{\Lambda}}  (\bar{T} , \bar{T})=o$, and this completes the proof of Part~(i).

Part~(ii): Let $\Gamma : = \End_\Lambda (T)^\op$ and regard $T$ as a $\Gamma^\op$-module in the usual way. 
Note that $R$ is noetherian by the hypothesis
and the $R$-algebra $\Gamma$ is finitely generated as a module over $R$. Therefore, $\Gamma$ is noetherian and the $\Gamma^\op$-module $T$ is finitely generated. In particular, $T$ is strongly finitely presented as a module over $\Lambda$ and $\Gamma^\op$. Furthermore, it is readily seen that $x_\bullet$ is regular on $\Gamma$ and the $\Gamma^\op$-module $T$. Now to prove the assertion it suffices to show that the $\Lambda$-module $T$ satisfies the conditions (i)--(iii) of Proposition~\ref{prop:tilting.characterization}. We prove this in three steps:

\textsc{Step~1.} Since $\pd_{\Lambda_x} (T_x) \leq n$ 
and $\pd_{\bar{\Lambda}} (T/xT) \leq n$ and $T$ 
is strongly finitely presented over $\Lambda$, $\pd_\Lambda (T) \leq n$ by Theorem~\ref{thm:rajabi}-(i).
Furthermore, Proposition~\ref{prop:self.orthogonal} shows that $\Ext^{\geqslant 1}_\Lambda (T,T)=o$, and so there are ring isomorphisms
$\Gamma / x_\bullet \Gamma \cong_\Gamma \End_{\bar{\Lambda}} (T/xT)^\op$ 
and $\Gamma_{x_\bullet} \cong \End_{\Lambda_x} (T_x)^\op$
by Lemma~\ref{lem:Ext.reg.element}-(ii) and  Lemma~\ref{lem:Ext.locx} respectively.

\textsc{Step~2.} By \textsc{Step~1}, $\Gamma_{x_\bullet}$ is tilted from $\Lambda_x$ via the
classical $n$-tilting $\Lambda_x$-module $T_x$ and 
$\Gamma / x_\bullet \Gamma$ is tilted from $\Lambda / x \Lambda$ via the classical $n$-tilting $\Lambda / x \Lambda$-module $T/xT$. By Proposition~\ref{prop:tilting.characterization}, 
\begin{equation} \tag{$\ast$}
\pd_{\Gamma^\op_{x_\bullet}} (T) \leq n \quad \text{and} \quad 
\pd_{\bar{\Gamma}^\op} (T/x_\bullet T) \leq n \: , 
\end{equation}
and
\begin{equation} \tag{$\ast \ast$}
\Ext^{\geqslant 1}_{\Gamma^\op_{x_\bullet}} (T_{x_\bullet} , T_{x_\bullet}) = \Ext^{\geqslant 1}_{\bar{\Gamma}^\op} (T/x_\bullet T , T/x_\bullet T)=o \: .
\end{equation}
Now ($\ast$) in conjunction with Theorem~\ref{thm:rajabi}-(i) implies that $\pd_{\Gamma^\op} (T) \leq n$, and ($\ast \ast$) in view of Lemma~\ref{lem:Ext.locx} implies that $\Ext^{\geqslant 1}_{\Gamma^\op} (T,T)=o$. 

\textsc{Step~3.} To complete the proof we need to show that the natural ring homomorphism $\rho: \Lambda \arrow \End_{\Gamma^\op} (T)$, 
which maps every $\lambda \in \Lambda$ to the $\Gamma^\op$-homomorphism $\lambda_\bullet : T \arrow T$ defined by multiplication by $\lambda$, is an isomorphism. 
Note that $\rho$ is also an $R$-homomorphism, which induces
the $R/xR$-homomorphism $\bar{\rho} : \Lambda / x \Lambda \arrow \End_{\Gamma^\op} (T) / x \End_{\Gamma^\op} (T)$ 
and $R_x$-homomorphism $\rho_x : \Lambda_x \arrow \End_{\Gamma^\op} (T)_x$.
By  Proposition~\ref{prop:tilting.characterization} we already know that the natural ring 
homomorphisms $\Lambda/ x \Lambda \stackrel{\cong}{\arrow} \End_{(\Gamma / x_\bullet \Gamma)^\op} (T/x_\bullet T)$ and
$\Lambda_x \stackrel{\cong}{\arrow} \End_{\Gamma^\op_{x_\bullet}} (T_{x_\bullet})$ are isomorphisms.
The ring isomorphism $\Lambda/ x \Lambda \stackrel{\cong}{\arrow} \End_{(\Gamma / x_\bullet \Gamma)^\op} (T/x_\bullet T)$ in view of 
Lemma~\ref{lem:Ext.reg.element}-(ii) implies that the $R/xR$-homomorphism $\bar{\rho}$ is an isomorphism and the ring isomorphism
$\Lambda_x \stackrel{\cong}{\arrow} \End_{\Gamma^\op_{x_\bullet}} (T_{x_\bullet})$ in view of Lemma~\ref{lem:Ext.locx} implies that
the $R_x$-homomorphism $\rho_x$ is an isomorphism. Therefore, 
\begin{align}
(\ker \rho)_x & = \ker (\rho_x) = o \label{eq:ker.of.rho} \\
(\coker \rho)_x & = \coker (\rho_x) = o \label{eq:coker.rho.x} \: ,
\end{align}
and
\begin{equation} \label{eq:cokerrho.mod.x}
\frac{\coker \rho}{x  (\coker \rho )} = \coker \bar{\rho} =o \:  .
\end{equation}
Since $\ker \rho$ is a submodule of $\Lambda$ and $x$ is non-zero-divisor on $\Lambda$,
$x$ is also a non-zero-divisor on $\ker \rho$ and therefore $\ker \rho = o$ in view of~\eqref{eq:ker.of.rho}.
Furthermore, Observation~\ref{obs:zerox} in conjunction with~\eqref{eq:coker.rho.x} and~\eqref{eq:cokerrho.mod.x} implies that $\coker (\rho)=o$.
Therefore, $\rho$ is a ring isomorphism.
\eproof

\brem \label{rem-der-equiv} 
Let $\Lambda$ be a ring and denote by $\D (\Lambda)$ the derived category of $\Lambda$.
It is well known~(see e.g.~\cite{rickard.der.morita} and~\cite{KZ}) that given a classical $n$-tilting triple $(\Lambda , T , \Gamma)$, there are mutually inverse equivalences
\[\xymatrix{\D (\Lambda) \ar@<0.6ex>[rr]^-{\Rhom_\Lambda (T,-)} \ar@{}[rr] |\sim & & \D (\Gamma) \ar@<0.6ex>[ll]^-{T \tensorL_{\Gamma} -}} \: . \]
Now if $x \in \Lambda$ is a central element which is regular on $\Lambda$ and $T$, then by Theorem~\ref{main.thm.classic} we have also the mutually inverse equivalences
\[\xymatrix{
\D (\Lambda_x) \ar@<0.6ex>[rr]^-{\Rhom_{\Lambda_x} (T_x,-)} \ar@{}[rr] |\sim & & \D (\Gamma_{x_\bullet}) \ar@<0.6ex>[ll]^-{T_x \tensorL_{\Gamma_{x_\bullet}} -} \\
\D (\bar{\Lambda}) \ar@<0.6ex>[rr]^-{\Rhom_{\bar{\Lambda}} (\bar{T},-)} \ar@{}[rr] |\sim & & \D (\bar{\Gamma}) \ar@<0.6ex>[ll]^-{\bar{T} \tensorL_{\bar{\Gamma}} -}
}\]
where $\bar{\Lambda}:=\Lambda/ x \Lambda$, $\bar{\Gamma}:=\Gamma/x_\bullet \Gamma$ and $\bar{T}:=T/xT$. 
This establishes a close relationship between various properties or invariants of $\Lambda_x$ and $\Gamma_{x_\bullet}$,
as well as $\bar{\Lambda}$ and $\bar{\Gamma}$. For example, it follows from the above equivalences or the Tilting Theorem (see~\cite[Corollary 2.4]{miyashita})
that
\begin{align*}
 \gldim (\Lambda) - \pd_\Lambda (T) & \leq \gldim (\Gamma) \leq \gldim (\Lambda) + \pd_\Lambda (T) \: , \\
\gldim (\Lambda_x) - \pd_{\Lambda_x} (T_x) & \leq \gldim (\Gamma_{x_\bullet}) \leq \gldim (\Lambda_x) + \pd_{\Lambda_x} (T_x)  \: , \\
\gldim (\Lambda / x \Lambda) - \pd_{\Lambda / x \Lambda} (T/xT) & \leq 
\gldim (\Gamma/x_{\bullet} \Gamma) \leq \gldim (\Lambda/x \Lambda) + \pd_{\Lambda/x \Lambda} (T/xT) \: .
\end{align*}
In particular, we have the following bi-implications:
\begin{align*}
\gldim (\Lambda) < + \infty & \iff \gldim (\Gamma) < + \infty \: ,  \\
\gldim (\Lambda_x) < + \infty & \iff \gldim (\Gamma_{x_\bullet}) < + \infty \: ,  \\
\gldim (\bar{\Lambda}) < + \infty & \iff 
\gldim (\bar{\Gamma}) < + \infty \: . 
\end{align*}
Furthermore, when $\Lambda$ is finitely generated over a noetherian center, the above inequalities in conjunction with Theorem~\ref{thm:rajabi} 
yield the following bi-implications:
\[\xymatrix{
\gldim (\Lambda) < + \infty \ar@{<=>}[r]  \ar@{<=>}[d] & 
\gldim (\Gamma) < + \infty \ar@{<=>}[d] \\
\gldim (\Lambda_x) ,  \gldim (\bar{\Lambda}) < + \infty \ar@{<=>}[r] & \gldim (\Gamma_{x_\bullet}) , \gldim (\bar{\Gamma}) < + \infty
}\]
\erem

We are now going to study non-classical tilting modules
under special base changes. As we shall see, the required
arguments in this case are more subtle due to lack of usual finiteness conditions on the tilting module under consideration, but this defect can be
compensated by the ``finite type property'' of tilting modules (recall~\ref{thm:finite.type}). We begin with the following 
non-classical version of Theorem~\ref{main.thm.classic}-(i): 

\begin{theorem} \label{main.thm.nonclassic.loc}
Let $\Lambda$ be a ring and $x \in \Lambda$ be a central element.
Let $T$ be a $\Lambda$-module on which $x$ is regular, and let $\Gamma:=\End_\Lambda (T)^\op$. Denote by $x_\bullet$ the 
central element of $\Gamma$ defined by $t \mapsto x t$ for every $t \in T$. If $(\Lambda , T , \Gamma)$ is an $n$-tilting triple, then $\big (\Lambda/x\Lambda , T/xT , \Gamma / x_\bullet \Gamma \big )$
and $(\Lambda_x , T_x , \Gamma_{x_\bullet})$ are $n$-tilting triples too.
\end{theorem}

\bproof
Let $\bar{\Lambda} : = \Lambda / x \Lambda$ and $\bar{T}:= T / x T$. Notice that since $T$ is an $n$-tilting module and $x$ is regular on $T$,
the property (T3) in the definition of a tilting module implies that $x$ is regular on $\Lambda$ too. This, in conjunction with flatness of $\Lambda_x$ over $\Lambda$,
implies that $\Tor^{\Lambda}_{\geqslant 1} (\Lambda_x , T) =\Tor^{\Lambda}_{\geqslant 1} (\bar{\Lambda} , T)=o$. 
Thus, in order to prove that $\big (\Lambda/x\Lambda , T/xT , \Gamma / x_\bullet \Gamma \big )$ and $(\Lambda_x , T_x , \Gamma_{x_\bullet})$ are $n$-tilting triples, it suffices by Proposition~\ref{prop:miyashita.base.change}, to show that for every cardinal number $\kappa$,
\[\Ext^{\geqslant 1}_{\bar{\Lambda}} \big (\bar{T} , \bar{T}^{(\kappa)} \big )=\Ext^{\geqslant 1}_{\Lambda_x} \big (T_x , (T_x)^{(\kappa)} \big )=o \: ,\]
and the canonical ring homomorphisms
$\Gamma / x_\bullet \Gamma \arrow \End_{\bar{\Lambda}} (\bar{T})$ and $\Gamma_{x_\bullet} \arrow \End_{\Lambda_x} (T_x)$ are isomorphisms.
We prove these in several steps.

\textsc{Step~1.} Since $\Ext^1_\Lambda (T,T)=0$ and
$x$ is regular on both $\Lambda$ and $T$, the canonical ring homomorphism $\Gamma / x_\bullet \Gamma \arrow \End_{\bar{\Lambda}} (\bar{T})$ is an isomorphism and the $\mathbb{Z}$-isomorphisms
\begin{equation} \label{eq:Ext.temp1}
\Ext^{n}_{\bar{\Lambda}} \big ( \bar{T} , \bar{T}^{(\kappa)} \big ) \cong \Ext^{n}_{\Lambda} \Big ( T , (T/xT)^{(\kappa)} \Big ) 
\cong \Ext^{n}_{\Lambda} \big ( T , \tfrac{T^{(\kappa)}}{x T^{(\kappa)}} \big )
\end{equation}
hold for every integer $n \geq 1$ by Lemma~\ref{lem:Ext.reg.element}.
On the other hand, the exact sequence
\[\xymatrix{ o \ar[r] & T \ar[r]^-x & T \ar[r] & \frac{T}{xT} \ar[r] & o }\]
of $\Lambda$-modules induces the exact sequence
\begin{equation} \label{eq:Ext.temp2}
\xymatrix{\Ext^n_\Lambda (T,T^{(\kappa)})  \ar[r] & \Ext^n_\Lambda \big ( T, \frac{T^{(\kappa)}}{x T^{(\kappa)}} \big ) \ar[r] & \Ext^{n+1}_\Lambda (T,T^{(\kappa)}) 
}\end{equation}
of $\mathbb{Z}$-modules for every integer $n \geq 1$ and every cardinal $\kappa$. Now~\eqref{eq:Ext.temp2} in view of $\Ext^{\geqslant 1}_\Lambda (T,T^{(\kappa)})=o$ and~\eqref{eq:Ext.temp1}
implies that $\Ext^{\geqslant 1}_{\bar{\Lambda}} \big ( \bar{T} , \bar{T}^{(\kappa)} \big )=o$. 

\textsc{Step~2.} Since the class $T^{\perp_\infty}$ is $n$-tilting, it follows from
Theorem~\ref{thm:finite.type} that there exists a resolving subcategory $\mathcal{S}$ of $\mod (\Lambda)$ contained in $\mathcal{P}_n^{< \omega}$ 
such that $T^{\perp_\infty} = \mathcal{S}^{\perp}$. If $\mathbf{\Omega}$ is the set of all syzygy modules of $T$ with respect to some projective
resolution of $T$, then the equality $T^{\perp_\infty} = \mathcal{S}^{\perp}$ implies ${}^\perp (\mathcal{S}^\perp) = {}^\perp (\mathbf{\Omega}^\perp )$ by Corollary~\ref{coro:T.perp.preenv}. From this one can deduce using Theorem~\ref{thm:eklof.trlifaj}-(ii) that:
\begin{description}
\item[(I)] Since $\Ext^{\geqslant 1}_\Lambda (T,T)=o$, $T \in {}^\perp (\mathbf{\Omega}^\perp )$ and so $T$ is a direct summand of some $\mathcal{S}$-filtered module $M$.

\item[(II)] Since $\mathcal{S} \subseteq {}^\perp (\mathcal{S}^\perp)={}^\perp (\mathbf{\Omega}^\perp )$, every element of $\mathcal{S}$ is a direct summand of some $\mathbf{\Omega}^+$-filtered module 
where $\mathbf{\Omega}^+ := \mathbf{\Omega} \cup \{ \Lambda  \}$.
\end{description}
We use these facts in the next two steps to show that the canonical
ring homomorphism $\Gamma_{x_\bullet} \arrow \End_{\Lambda_x} (T_x)$ is an isomorphism, and that $\Ext^{\geqslant 1}_{\Lambda_x} \big (T_x , (T_x)^{(\kappa)} \big )=o$ for every cardinal $\kappa$.

\textsc{Step~3.} In order to prove that the canonical homomorphism $\Gamma_{x_\bullet} \arrow \End_{\Lambda_x} (T_x)$ is an isomorphism,
it suffices by part~(I) of \textsc{Step~2} to show  that the canonical $\mathbb{Z}$-homomorphism
$\theta_M : \hom_\Lambda (M,T)_x \arrow \hom_{\Lambda_x} (M_x , T_x)$ is an isomorphism.
Let $\{M_\alpha \}_{\alpha \leq \sigma}$ be a $\mathcal{S}$-filtration of the $\Lambda$-module $M$.
In order to prove $\theta_M$ is a $\mathbb{Z}$-isomorphism, we may prove by transfinite induction that for every ordinal number $\alpha$ with $\alpha \leq \sigma$, the canonical $\mathbb{Z}$-homomorphism $\theta_{M_\alpha} :\hom_\Lambda (M_\alpha ,T)_x \arrow \hom_{\Lambda_x} \big ( (M_\alpha)_x , T_x \big ) $ is an isomorphism.
Since $M_0=o$, it is clear that $\theta_{M_\alpha}$ is a $\mathbb{Z}$-isomorphism for $\alpha=0$. Assume now that $\theta_{M_\alpha}$ is an
isomorphism for some ordinal $\alpha < \sigma$, and consider the short exact sequence
\begin{equation} \label{eq:A1}
\xymatrix{o \ar[r] & M_\alpha \ar[r] & M_{\alpha + 1} \ar[r] & M_{\alpha+1} / M_{\alpha} \ar[r] & o }
\end{equation}
of $\Lambda$-modules as well as the induced short exact sequence
\begin{equation} \label{eq:A2}
\xymatrix{o \ar[r] & (M_\alpha)_x \ar[r] & (M_{\alpha + 1})_x \ar[r] & \big ( M_{\alpha+1} / M_{\alpha} \big )_x \ar[r] & o }
\end{equation}
of $\Lambda_x$-modules. Applying the functor $\hom_\Lambda (-,T)$ to~\eqref{eq:A1} and
the functor $\hom_{\Lambda_x} (-,T_x)$ to~\eqref{eq:A2} we obtain the
exact sequences
\begin{equation} \label{eq:A3}
\begin{scriptsize}
\xymatrixcolsep{0.5cm}
\xymatrix{o \ar[r] & \hom_\Lambda \big ( \frac{M_{\alpha+1}}{M_\alpha} , T \big ) \ar[r] & \hom_\Lambda (M_{\alpha+1} , T) \ar[r] &
\hom_\Lambda (M_\alpha , T) \ar[r] & \Ext^1_\Lambda \big ( \frac{M_{\alpha+1}}{M_\alpha} , T \big )}
\end{scriptsize}
\end{equation}
and
\begin{equation}\label{eq:A4}
\begin{scriptsize}
\begin{array}{c}
\xymatrixcolsep{0.6cm}
\xymatrix{o \ar[r] & \hom_{\Lambda_x}  \big ( (\frac{M_{\alpha+1}}{M_\alpha})_x , T_x \big ) \ar[r] &
\hom_{\Lambda_x}  \big ( (M_{\alpha+1})_x , T_x \big ) \ar[r] &
\hom_{\Lambda_x} \big ( (M_\alpha)_x , T_x \big )   \\
 &  & \ar[r] & \Ext^1_\Lambda \big ( (\frac{M_{\alpha+1}}{M_\alpha})_x , T_x \big )} 
\end{array}
\end{scriptsize}
\end{equation}
of $\Lambda$-modules. The $\Lambda$-module $M_{\alpha+1}/ M_{\alpha}$ is a direct summand of a $\mathbf{\Omega}^+_T$-filtered module by part~(II) in \textsc{Step~2}. On the other hand, every element of $\mathbf{\Omega}^+_T$ belongs to ${}^\perp T$ because $\Ext^{\geqslant 1}_\Lambda (T,T)=o$.
Therefore, $M_{\alpha+1}/ M_{\alpha}$ belongs to ${}^\perp T$ by Eklof Lemma~\ref{Eklof.lemma}. Thus
$\Ext^1_\Lambda \big ( \frac{M_{\alpha+1}}{M_\alpha} , T \big )=o$, and from this one concludes in view of Lemma~\ref{lem:Ext.locx} that
\[ \Ext^1_{\Lambda} \big ( M_{\alpha+1} / M_\alpha , T \big )_x \cong_{\mathbb{Z}} 
\Ext^1_{\Lambda_x} \big ( (M_{\alpha+1} / M_\alpha)_x , T_x \big ) = o \: . \]
Therefore, the exact sequences~\eqref{eq:A3} and~\eqref{eq:A4} yield the commutative diagram
\begin{small}
\[
\xymatrixcolsep{0.5cm}
\xymatrix{
o \ar[r] & \hom_\Lambda \big ( \frac{M_{\alpha+1}}{M_\alpha} , T \big )_x \ar[r] \ar[d]_-{\theta_{\frac{M_{\alpha+1}}{M_\alpha}}} & \hom_\Lambda (M_{\alpha+1} , T)_x \ar[r] \ar[d]_-{\theta_{M_\alpha+1}} &
 \hom_\Lambda(M_\alpha , T)_x \ar[r] \ar[d]_-{\theta_{M_{\alpha}}} & o \\
o \ar[r] &  \hom_{\Lambda_x} \big ( (\frac{M_{\alpha+1}}{M_\alpha})_x , T_x \big ) \ar[r] & \hom_{\Lambda_x} \big ( (M_{\alpha+1})_x , T_x \big ) \ar[r] &
\hom_{\Lambda_x} \big ( (M_\alpha)_x , T_x \big ) \ar[r] & o
}\]
\end{small}
of $\mathbb{Z}$-modules with exact rows. The map $\theta_{M_\alpha}$ is an isomorphism by the inductive hypothesis and the map $\theta_{\frac{M_{\alpha+1}}{M_\alpha}}$ is an isomorphism by Lemma~\ref{lem:Ext.locx}. Therefore,
$\theta_{M_{\alpha+1}}$ is an isomorphism by Short Five Lemma. Finally, if $\alpha \leq \sigma$ 
is a limit ordinal and $\theta_{M_\beta}$ is an isomorphism for all $\beta < \alpha$, 
then $M_\alpha = \varinjlim_{\beta < \alpha} M_\beta$ which implies that
$(M_\alpha)_x = \varinjlim_{\beta < \alpha} (M_\beta)_x$. Therefore, we have the following
commutative diagram
\[\xymatrix{
\hom_\Lambda (M_\alpha , T)_x \ar[rr]^-{\theta_{M_\alpha}} \ar[d]^-{\cong} && \hom_{\Lambda_x} \big ( (M_\alpha)_x , T_x \big ) \ar[d]^-{\cong} \\
\varprojlim \limits_{\beta< \alpha} \hom_{\Lambda} (M_\beta , T)_x \ar[rr]^-{\varprojlim \limits_{\beta < \alpha} \theta_{M_\beta}} && 
\varprojlim \limits_{\beta < \alpha} \hom_{\Lambda_x} \big ( (M_\beta)_x , T_x \big ).
}\]
of $\mathbb{Z}$-modules where the vertical arrows are the canonical isomorphisms. Since by the hypothesis each $\theta_{M_\beta}$ is an isomorphism for every $\beta < \alpha$, the map $\varprojlim_{\beta < \alpha} \theta_{M_\beta}$ and hence $\theta_{M_\alpha}$ is an isomorphism. This implies by transfinite induction that the canonical map $\theta_M$ is an isomorphism, and this completes the proof of the fact that
the canonical ring homomorphism $\Gamma_{x_\bullet} \arrow \End_{\Lambda_x} (T_x)$ is an isomorphism.

\textsc{Step~4.} We conclude the proof by showing that $\Ext^{\geqslant 1}_{\Lambda_x} \big (T_x , T_x^{(\kappa)} \big )=o$
for every cardinal number $\kappa$. By part~(I) in \textsc{Step~2}, the $\Lambda$-module $T$ is a direct summand of a $\mathcal{S}$-filtered module $M$.
Since the $\Lambda_x$-module $T_x$ is a direct summand of $M_x$, we may prove $\Ext^{\geqslant 1}_{\Lambda_x} \big (T_x , T_x^{(\kappa)} \big )=o$ by showing that $\Ext^{\geqslant 1}_{\Lambda_x} \big (M_x , T_x^{(\kappa)} \big )=o$ or equivalently $M_x \in {}^\perp \big (T_x ^{(\kappa)}\big )$.
Let $\mathcal{S}_x : = \{ C_x : C \in \mathcal{S} \}$. Since $M$ is $\mathcal{S}$-filtered and $\Lambda_x$ is flat over $\Lambda$,
the $\Lambda_x$-module $M_x$ is $\mathcal{S}_x$-filtered. 
Note that every element of $\mathcal{S}$ is a direct summand of 
a $\mathbf{\Omega}^+$-filtered module and every element of $\mathbf{\Omega}^+$ belongs  to ${}^\perp (T^{(\kappa)})$ because $\Ext^{\geqslant 1}_{\Lambda} (T,T^{(\kappa)})=o$.
Therefore, $\mathcal{S} \subseteq {}^\perp (T^{(\kappa)})$ by Eklof Lemma~\ref{Eklof.lemma}, which implies by
Lemma~\ref{lem:Ext.locx} that $\mathcal{S}_x \subseteq {}^\perp (T_x^{(\kappa)})$. Therefore, the $\Lambda_x$-module $M_x$ is ${}^\perp (T_x^{(\kappa)} )$-filtered and hence belongs to ${}^\perp (T_x^{(\kappa)})$ by Eklof Lemma~\ref{Eklof.lemma}.
The proof is thus complete.
\eproof

Unlike the classical case (cf. Theorem~\ref{main.thm.classic}-(ii)), the converse of Theorem~\ref{main.thm.nonclassic.loc} does \emph{not} hold in general,
as the following example shows.

\bex \label{counterexample}
Let $\Lambda:=\mathbb{Z}_{(2)}$ be the localization of $\mathbb{Z}$ at the prime ideal $(2)$ and note that $\mathbb{Z}_{(2)}$ is a DVR whose
field of fractions is $\mathbb{Q}$. Let $T:=\mathbb{Q} \amalg \mathbb{Z}_{(2)}$ and note that the central element
$x:=2 \in \Lambda$ is regular on both $\Lambda$ and $T$.
Since $\Lambda / x \Lambda$ and $\Lambda_x$ are fields, the $\Lambda / x \Lambda$-module $T/xT$ and the $\Lambda_x$-module $T_x$
are tilting while the $\Lambda$-module $T$ is \emph{not} self-orthogonal and hence \emph{cannot} be a (partial) tilting module.
\eex

However, we can prove the following partial converse to Theorem~\ref{main.thm.nonclassic.loc} under some relatively reasonable conditions.

\bthm \label{main.thm:tilting.class}
Let $\Lambda$ be a ring over which every $\Lambda$-module of finite flat dimension is also of finite projective dimension (see Remark~\ref{rem:inf.perf}). Let  $T$ be a $\Lambda$-module such that $T^{\perp_\infty}$ is closed under
direct sums and let $x \in \Lambda$ be a central element which is regular on $\Lambda$ and $T$. 
\begin{enumerate}
\item If the $\Lambda_x$-module $T_x$ 
and the $\Lambda / x \Lambda$-module $T/xT$ are of finite projective dimension, then $T^{\perp_\infty}$ is a tilting class.

\item Assume that $\End_\Lambda (T)$ is left-noetherian and
$T$ is finitely generated as an $\End_\Lambda (T)$-module.
If the $\Lambda_x$-module $T_x$ and the $\Lambda / x \Lambda$-module $T/xT$ are tilting, then $T$ is a partial tilting module. In particular, $T$ is a direct summand of a tilting module.
\end{enumerate}
\ethm

\bproof
Let $\bar{\Lambda} := \Lambda / x \Lambda$ and $\bar{T} := T / x T$.

Part~(i): Since by the hypothesis $\pd_{\Lambda_x} (T_x)$ and $\pd_{\bar{\Lambda}} (\bar{T})$ are finite, it follows from Theorem~\ref{thm:rajabi}-(i) and the hypothesis
that $n:=\pd_\Lambda (T) < +\infty$. In particular, $T \in \mathcal{P}_n$ implies ${}^\perp (T^{\perp_\infty}) \subseteq {}^\perp (\mathcal{P}^\perp_n)=\mathcal{P}_n$. The class $T^{\perp_\infty}$ is closed under direct sums by the hypothesis, and it is also coresolving and closed under direct summands. Furthermore, $T^{\perp_\infty}$ is special preenveloping by Corollary~\ref{coro:T.perp.preenv}. Consequently, $T^{\perp_\infty}$ is an $n$-tilting class by Theorem~\ref{thm:tilting.class}.

Part~(ii): The proof of this part proceeds in two steps:

\textsc{Step~1.} Since the $\Lambda_x$-module $T_x$ 
and the $\Lambda / x \Lambda$-module $T/xT$ are tilting, 
the class $T^{\perp_\infty}$ is tilting by part~(i) of the theorem.
Thus, it follows from Theorem~\ref{thm:finite.type} that $T^{\perp_\infty} = \mathcal{S}^\perp$ for some resolving subcategory $\mathcal{S}$ of $\mod (\Lambda)$ contained in $\mathcal{P}^{< \omega}_n$.
Let $\mathbf{\Omega}_T$ be the set of syzygy modules of $T$ with respect to a projective resolution $\mathbf{P}_T$ of $T$,
and note that ${}^\perp \big ( (\mathbf{\Omega}_T)^\perp \big ) = {}^\perp \big ( \mathcal{S}^\perp \big )$ by Corollary~\ref{coro:T.perp.preenv}.
Therefore, by Theorem~\ref{thm:eklof.trlifaj}~(ii),
every element of $\mathcal{S}$ is a direct summand of a $\mathbf{\Omega}^+_T$-filtered module, where $\mathbf{\Omega}^+_T:=\mathbf{\Omega}_T \cup \{ \Lambda  \}$.
Note also that $\Tor^{\Lambda}_{\geqslant 1} (\Lambda_x , T) = o$ (because $\Lambda_x$ is flat over $\Lambda$) and 
$\Tor^{\Lambda}_{\geqslant 1} (\bar{\Lambda} , T) =o$ (because $x$ is regular on $\Lambda$ and $T$). Therefore,
$\mathbf{P}_{T_x} : = \Lambda_x \otimes_\Lambda \mathbf{P}$ is a projective resolution of the $\Lambda_x$-module $T_x$ and $\mathbf{P}_{\bar{T}}:=\bar{\Lambda} \otimes_\Lambda \mathbf{P}$ is a projective resolution of the $\bar{\Lambda}$-module $\bar{T}$. Let $\mathbf{\Omega}_{T_x}$ be the set of syzygy modules 
of the $\Lambda_x$-module $T_x$ with respect to $\mathbf{P}_{T_x}$, and likewise let $\mathbf{\Omega}_{\bar{T}}$ be 
the set of syzygy modules of the $\bar{\Lambda}$-module $\bar{T}$ with respect to $\mathbf{P}_{\bar{T}}$.

\textsc{Step~2.} By the preceding step we already know that the $\Lambda$-module $T$ satisfies conditions (P1) and (P3) in the definition of a partial tilting module; cf.~\ref{df:partial.tilting}. Thus, in order to complete the proof, we only need to show that for every cardinal number $\kappa$, $\Ext^{\geqslant 1}_\Lambda (T,T^{(\kappa)})=o$ or equivalently $T^{(\kappa)} \in T^{\perp_\infty}$.
But this amounts to $\Ext^1_\Lambda (C , T^{(\kappa)})=o$ for every $C \in \mathcal{S}$, because $T^{\perp_\infty}=\mathcal{S}^\perp$.
Let $C \in \mathcal{S}$ and note that $C$ is a direct summand of a $\mathbf{\Omega}^+_T$-filtered $\Lambda$-module.
Since $\Lambda_x$ is flat over $\Lambda$,  the $\Lambda_x$-module $C_x$ is a direct summand of a $\mathbf{\Omega}^+_{T_x}$-filtered module,
where $\mathbf{\Omega}_{T_x}^+:=\mathbf{\Omega}_{T_x} \cup \{ \Lambda_x \}$. Furthermore, since $x$ is non-zero-divisor on every member of $\mathbf{\Omega}^+_T$, an easy argument based on transfinite induction shows that $x$ is also non-zero-divisor on $C$. Thus
$\Tor^{\Lambda}_{\geqslant 1} (\bar{\Lambda} , C)=o$ which implies that the $\bar{\Lambda}$-module $\bar{C}$
is a direct summand of a $\mathbf{\Omega}^+_{\bar{T}}$-filtered module, where 
$\mathbf{\Omega}_{\bar{T}}^+:=\mathbf{\Omega}_{\bar{T}} \cup \{ \bar{\Lambda} \}$. On the other hand, since by the hypothesis $T_x$ and $\bar{T}$ are tilting modules, 
$\Ext^{\geqslant 1}_{\Lambda_x} \big (  T_x , T_x^{(\kappa)} \big )=\Ext^{\geqslant 1}_{\bar{\Lambda}} \big ( \bar{T} , \bar{T}^{(\kappa)}  \big )=o$.
Therefore, $\mathbf{\Omega}_{T_x} \subseteq {}^{\perp} \big ( T_x^{(\kappa)} \big )$ 
and $\mathbf{\Omega}_{\bar{T}} \subseteq {}^{\perp} (\bar{T}^{(\kappa)})$. These inclusions in view of Eklof Lemma~\ref{Eklof.lemma}, imply that
$C_x \in {}^{\perp} \big ( T_x^{(\kappa)} \big )$ and $\bar{C} \in {}^{\perp} (\bar{T}^{(\kappa)})$.
That is, $\Ext^{1}_{\Lambda_x} (C_x , (T_x)^{(\kappa)})=\Ext^{1}_{\bar{\Lambda}} \big ( \bar{C} , (\bar{T})^{(\kappa)} \big )=o$ which implies by Proposition~\ref{prop:self.orthogonal} that $\Ext^1_\Lambda (C , T^{(\kappa)})=o$.

We have thus proved that the $\Lambda$-module $T$ is a partial tilting module. In particular, $T$ is a direct summand of some tilting module by Remark~\ref{rem:partial}.
\eproof

\brem \label{rem:inf.perf}
A couple of remarks are in order regarding the hypotheses of Theorem~\ref{main.thm:tilting.class}:
\begin{enumerate}
\item 
Rings over which every module of finite flat dimension is also of finite projective dimension are abundant. Examples are: perfect rings (in particular, artin algebras), rings with finite (big) finitistic projective dimension (in particular, commutative noetherian rings of finite Krull dimension) by~\cite[Proposition~6]{jensen.Lim}, and right-noetherian algebras admitting a ``dualizing complex'' (see~\cite{jorgensen.finiteflat}).
\item 
The conditions imposed on $T$ and $\End_{\Lambda}(T)$ seem to be \emph{relatively} natural:
Indeed, ``being finitely generated over the endomorphism ring" is a necessary condition for a module to be tilting (this could easily be deduced from \cite[Proposition~(1.2)]{tilting.preenv}). This, to some extent, justifies the finiteness condition imposed on $T$. 
On the other hand, the endomorphism ring of a (partial) tilting module is \emph{not} necessarily left-noetherian (consider a non-right-noetherian ring as a tilting module). To the best of our knowledge, little is known about chain conditions of endomorphism rings of non-classical tilting or partial tilting modules. Thus, if nothing else, one can give our hypotheses the benefit of the doubt.
%
\end{enumerate} 
\erem

We conclude with the following application of 
our results, which illustrates a kind of ``tilting exchange" among rings.

\bex
Suppose that $\Lambda$ is a filtered ring with the filtration $\cdots \subset F_n\Lambda \subset F_{n+1}\Lambda \subset \cdots$.
Let $\mathfrak{R}(\Lambda):=\bigoplus_{n \in \mathbb{Z}} F_n\Lambda$ and $\operatorname{Gr }(\Lambda) := \bigoplus_{n \in \mathbb{Z}} \frac{F_n\Lambda}{F_{n-1}\Lambda}$
be the \emph{Rees ring} and
the \emph{associated graded ring} of $\Lambda$ respectively.
The element $1 \in F_1\Lambda$ corresponds to a central regular element $x$ in $\mathfrak{R}(\Lambda)$, and there are ring isomorphisms
$\mathfrak{R}(\Lambda)_x \cong \Lambda[x,x^{-1}]$ 
and $\mathfrak{R}(\Lambda)/x\mathfrak{R}(\Lambda) \cong \operatorname{Gr }(\Lambda)$.
Now suppose that $T$ is a graded module over $\mathfrak{R}(\Lambda)$ and that $x$ is regular on $T$.
Then by our results (cf. Theorem~\ref{main.thm.classic} and Theorem~\ref{main.thm.nonclassic.loc}), if $T$ is a (classical) tilting module over $\mathfrak{R}(\Lambda)$, then $T_x$ is a (classical) tilting module over $\Lambda[x,x^{-1}]$ and $T/xT$ is a (classical) tilting module over $\operatorname{Gr}(\Lambda)$. 
It is worth noting that under our assumption, the tilting $\mathfrak{R} (\Lambda)$-module $T$ is ``extended'' in the sense that there is a filtered $\Lambda$-module $M$ which is tilting and 
its associated Rees module $\mathfrak{R}(M)$ is isomorphic to $T$.
Indeed, the functor $\mathfrak{R}(-)$ from the category of filtered $\Lambda$-modules to the category of graded 
$\mathfrak{R} (\Lambda)$-modules is fully faithful and its essential image is the full subcategory of $x$-torsion-free graded modules (see e.g. \cite{Bergh}). As $T$ lies in the essential image of $\mathfrak{R} (-)$, there exists a filtered $\Lambda$-module $M$ with $\mathfrak{R} (M) \cong T$. On the other hand, there is a ring isomorphism $\Lambda \cong \mathfrak{R}(\Lambda)/(1-x)\mathfrak{R}(\Lambda)$ under which
$M \cong \mathfrak{R}(M)/(1-x)\mathfrak{R}(M) \cong T/(1-x)T$ as $\Lambda$-modules (c.f. \emph{loc. cit.}). 
Now, considering the grading, it is clear that $1-x$ is a central element which is regular on both $\mathfrak{R}(\Lambda)$ and $T$. Thus Theorem~\ref{main.thm.nonclassic.loc} shows that $M$ is a tilting module over $\Lambda$.
\eex

\section*{Acknowledgment}

The authors would like to thank the anonymous referee for reading the manuscript carefully and spotting a gap in the proof of Theorem~\ref{main.thm:tilting.class} in an earlier version of the manuscript. The research of the second author is supported by Iran National Science Foundation (INSF), Project No. 94025823. The research of the third author is supported by Iran National Science Foundation (INSF), Project No.~95831492.

\bibliographystyle{amsplain} 
\bibliography{biblio}
\end{document}